\documentclass{amsart}
\usepackage{amssymb}
\usepackage{graphicx}
\def\({\left(}
\def\){\right)}

\newtheorem{lema}{Lemma}[section]

\newtheorem*{teorema*}{Theorem}
\newtheorem{proposicion}[lema]{Proposition}

\newtheorem{example}[lema]{Example}

\newtheorem{corollary}[lema]{Corollary}
\newtheorem{theorem}[lema]{Theorem}

\newtheorem{proposition}[lema]{Proposition}

\newtheorem{definition}[lema]{Definition}

  {\par \hfill \fbox{}}
  {\par \hfill \fbox{}}
  {\par \hfill }
  {\par \hfill }

\def\D{\mathcal D}
\def\DD{\mathbb D}

\def\NN{\mathbb N}
\def\RR{\mathbb R}
\def\TT{\mathbb T}

\def\beq{\begin{equation}}
\def\eeq{\end{equation}}
\def\beginpf{\noindent{\bf Proof.} \quad}
\def\endpf{\qquad\hfill\rule{2.2mm}{2.2mm}\break}

\def\ds{\displaystyle}
\def\epsilon{\varepsilon}
\newcommand\clos{\mathop{\rm clos}\nolimits}
\newcommand\supp{\mathop{\rm supp}\nolimits}
\begin{document}

\title[Weighted composition operators]{Boundedness, compactness and Schatten-class \\
membership of weighted composition operators}
\author{Eva A. Gallardo-Guti\'errez}
\address{Departamento de Matem\'aticas, Universidad de Zaragoza e IUMA,
Plaza San Francisco s/n, 50009 Zaragoza, SPAIN.}
\email{eva@unizar.es}
\author{Romesh Kumar}
\address{Department of Mathematics, University of Jammu,
Jammu 180 001, INDIA.} \email{romesh\_jammu@yahoo.com}
\author{Jonathan R.
Partington}
\address{School of
Mathematics, University of Leeds, Leeds LS2 9JT, U.K.}
\email{J.R.Partington@leeds.ac.uk}

\thanks{The first and third authors are partially supported by Plan Nacional I+D grant no.
MTM2007-61446 and Go\-bier\-no de Arag\'on research group
\emph{An\'alisis Matem\'atico y Aplicaciones}, ref. DGA E-64.
The second author is partially supported by the Royal Society (UK) and the Department of Science and Technology (India).}

\date{July 2009}

\begin{abstract}
The boundedness and compactness of weighted composition operators
on the Hardy space ${\mathcal H}^2$ of the unit disc is analysed.
Particular reference is made to the case when the self-map of the
disc is an inner function. Schatten-class membership is also
considered; as a result, stronger forms of the two main results of
a recent paper of Gunatillake are derived. Finally, weighted
composition operators on weighted Bergman spaces
$\mathcal{A}^2_\alpha(\mathbb{D})$ are considered, and the results
of Harper and Smith, linking their properties to those of Carleson
embeddings, are extended to this situation.
\end{abstract}

\maketitle

\section{Introduction}

Let $\mathbb {D}$ denote the open unit disk of the complex plane.
For $1\leq p< \infty$, recall that the classical Hardy space
$\mathcal{H}^p$ consists of holomorphic functions $f$ on $\mathbb
D$ for which the norm
$$
\|f\|_p=\left (\sup_{0\leq r<1} \int_{0}^{2\pi}
|f(re^{i\theta})|^p \, \frac{d\theta}{2\pi} \right )^{1/p}
$$
is finite. If $p=\infty$, $\mathcal{H}^{\infty}$ is the space of
holomorphic functions $f$ on $\mathbb{D}$ such that
$$
\|f\|_{\infty}=\sup_{\mathbb{D}} |f(z)|<\infty.
$$
Fatou's Theorem asserts that any Hardy function $f$ has  radial
limit at $e^{i\theta}\in \partial \mathbb{D}$ except on a set
Lebesgue measure zero (see \cite{Du}, for instance). Throughout
this work, $f(e^{i\theta})$ will denote the radial limit of $f$ at
$e^{i\theta}$, i. e., $f(e^{i\theta})=\lim_{r\to 1^{-}}
f(re^{i\theta})$.

If $\varphi$ is an analytic function on $\mathbb D$ which takes
$\mathbb{D}$ into itself, the Littlewood Subordination Principle
\cite{littlewood} ensures that the composition operator induced by
$\varphi$,
$$
C_{\varphi} f =f\circ \varphi,\qquad (f\in \mathcal{H}^p),
$$
is bounded on $\mathcal{H}^p$, $1\leq p \leq \infty$.

On the other hand, given $h \in \mathcal{H}^p$ it is possible to
define for those functions $f\in \mathcal{H}^p$ for which it makes
sense the weighted composition operator $W_{h,\varphi}$
\[
W_{h,\varphi}f (z) = h(z) f(\varphi(z)).
\]
Note that if $f$ belongs to the domain of the operator, which will
denoted by  $\D(W_{h,\varphi})$, the weighted composition operator
may be expressed by  $$W_{h,\varphi}=T_h C_\varphi,$$ where $T_h$
denotes the Toeplitz operator with symbol $h$. Observe that the
condition $h \in \mathcal{H}^\infty$ is always a sufficient
condition for boundedness of $W_{h,\varphi}$. By considering the
image of the constant functions, it is clear that $h \in
\mathcal{H}^p$ is a necessary condition. Nevertheless, the most
interesting feature in this sense is that there exists
symbols $h$ (unbounded) and $\varphi\in \mathcal{H}^{\infty}$ with
$\varphi(\mathbb{D})\subseteq \mathbb{D}$ such that
$W_{h,\varphi}$ is a bounded operator in $\mathcal{H}^p$.

Similar considerations, with the obvious modifications, apply to
other spaces of analytic functions such as the Bergman space
$\mathcal{A}^2(\mathbb{D})$. For some general properties of these
operators on Hardy and Bergman spaces we refer to the survey
\cite{KP}, where it is remarked that characterising the
boundedness of weighted composition operators on spaces of
functions on the disc enables one to characterise the boundedness
of such operators on Hardy and Bergman spaces of other
simply-connected domains, such as the half-plane.

Though weighted composition operators have attracted recently the
attention of operator theorists (see for instance the recent
papers \cite{CH}, \cite{CH1}, \cite{zen}, \cite{zenmartin}), we
would like to point out that its study traces back to the sixties.
Indeed, de Leeuw showed that the isometries in the Hardy space
$\mathcal{H}^1$ are weighted composition operators and Forelli
obtained the same result for the Hardy spaces $\mathcal{H}^p$ when
$1 < p < \infty$, $p\neq 2$ (see \cite{Ho} and \cite{Fo}).

The aim of this work is taking further the study of weighted
composition on spaces of analytic functions of $\mathbb{D}$. In
particular, we deal with the boundedness of weighted composition
operators on Hardy and Bergman spaces of the unit disc. In this
sense, we are primarily interested in the behaviour of the weights
$h$. Let us point out in this direction that an applicable
characterization of those weights $h$ so that $W_{h,\varphi}$ is
bounded on the Bergman space would solve the famous Brennan
conjecture (see Shimorin's work \cite{Shi} and the related work by
Smith \cite{Sm}).

In Section~\ref{sec:2} we use recent results of Harper and
Smith~\cite{zen,zenmartin} (of which the main tool is the
reproducing kernel of the Hilbert space $\mathcal{H}^2$) to derive
new results on the boundedness, compactness and Schatten-class
membership of $W_{h,\varphi}$. A case of particular interest is
when $\varphi$ is inner, and in this case the properties of the
angular derivative are seen to play an important role. In
particular, we give extensions of recent results due to
Jury~\cite{jury} and Gunatillake~\cite{gajath}. Most of our
results apply principally to operators on the Hilbert space
$\mathcal{H}^2$, but some results and proofs make sense also for
all $\mathcal{H}^p$ with $1 \le p < \infty$.

In Section~\ref{sec:3} we  extend some results concerning the link between weighted
composition operators and Carleson embeddings, due to Harper and Smith \cite{zen,zenmartin},
to the case of weighted Bergman spaces.

\section{Boundedness and compactness on Hardy spaces}
\label{sec:2}
\subsection{Conditions for boundedness}

Contreras and  Hern\'andez-D\'\i az \cite{CH} gave a necessary and
sufficient condition for boundedness of $W_{h,\varphi}$, namely
that $\mu_{h,\varphi}$ is a Carleson measure on $\overline\DD$,
where \beq\label{eq:carl} \mu_{h,\varphi}(E)=\int_{\varphi^{-1}(E)
\cap \TT} |h|^2 dm, \eeq for measurable subsets $E \subseteq \DD$.

Zen Harper \cite[Thm.~3.3]{zen} exploited the theory of
admissibility for semigroups to show (as part of a family of more
general results) that a more explicit necessary and sufficient
condition for boundedness of $W_{h,\varphi}$ is the following:
\beq\label{eq:zen} \sup_{|w|<1} \left \|
\frac{(1-|w|^2)^{1/2}h}{1-\overline w \varphi} \right \|_2 <
\infty. \eeq Moreover, $W_{h,\varphi}$ is compact if and only if
\beq\label{eq:zencompact} \left \|
\frac{(1-|w|^2)^{1/2}h}{1-\overline w \varphi} \right \|_2 \to 0
\qquad \hbox{as} \quad |w| \to 1. \eeq

Some independent generalizations to mappings between
$\mathcal{H}^p$ and $\mathcal{H}^q$  spaces were
given by Cu\u ckovi\'c and Zhao \cite{CZ}.

In fact, (\ref{eq:zen}) follows directly from
(\ref{eq:carl}) on using the standard reproducing kernel test for Carleson
embeddings (see \cite[p.~231]{garnett} or \cite[p.~105]{nik}), namely that an embedding $J:
\mathcal{H}^2 \to L^2(\mu)$ is bounded (i.e., $\mu$ is a Carleson
measure on $\DD$) if and only if $J$ is uniformly bounded on the
set of normalized reproducing kernels $\{\widetilde k_w: \, w \in
\DD\}$ on $\mathcal{H}^2$, where
\[
\widetilde k_w(z)=\frac{(1-|w|^2)^{1/2}}{1-\overline w z}, \qquad
\hbox{for} \quad z \in \DD.
\]

For
\[
\int_{\DD} |\widetilde k_w|^2 \, d\mu_{h,\varphi} =  \int_{\TT}
|h|^2 |\widetilde k_w \circ \varphi|^2 \, dm,
\]
from which (\ref{eq:zen}) follows easily. Similar arguments
involving vanishing Carleson measures establish
(\ref{eq:zencompact}).

On the other hand,  Jury \cite{jury} has recently obtained the
estimate $\|W_{h,\varphi}\| \le \|h\|_{H(\varphi)}$, but only for
$h \in H(\varphi)$, where $H(\varphi)$ is the de Branges--Rovnyak
space, that is, the reproducing kernel Hilbert space on $\DD$ with
kernel
\[
k^\varphi_w(z) =
\frac{1-\overline{\varphi(w)}\varphi(z)}{1-\overline w z}.
\]
In general this is not a necessary condition   for boundedness. If
$\varphi$ is inner, then $H(\varphi)=K_\varphi:= \mathcal{H}^2
\ominus \varphi \mathcal{H}^2$, and $P_{K_\varphi}h= \varphi
P_-(\overline\varphi h)$, where $P_-$ is the orthogonal projection
onto $L^2 \ominus \mathcal{H}^2$. In this case, it is easy to
check that, if $h \in K_\varphi$,
\[
\left \| \sum_{n=0}^\infty a_n \varphi^n h \right \|^2_2 =
\|h\|_2^2\sum_{n=0}^\infty |a_n|^2,
\]
since
\[
\langle \varphi^n h, \varphi^m h\rangle = \left
\{\begin{array}{ll} 0 &
\mbox{ for } n \ne m, \\
\noalign{\medskip}
 \|h\|_2^2 & \mbox{ for } n=m, \end{array} \right.
\]
and so $\| W_{h,\varphi} f\|_2=\|h\|_2 \|f\|_2$. Thus the
condition is correct for $h \in K_\varphi$, although not in
general. Indeed, if $\varphi(z)=z$, then
$\|W_{h,\varphi}\|=\|T_h\|=\|h\|_\infty$.

Suppose that $\varphi$ is inner, and let $\alpha$ be an
automorphism of the disc such that $\alpha(\varphi(0))=0$. Then
$\alpha \circ \varphi$ is also inner, and $W_{h,\varphi}C_\alpha=
W_{h,\alpha \circ \varphi}$. For questions of boundedness and
compactness, we may therefore assume without loss of generality
that $\varphi(0)=0$, and hence $\{1,\varphi,\varphi^2,\ldots \}$
is an orthonormal set and $C_\varphi$ is an isometry.

A further observation is that $W_{h,\varphi}$ is bounded as soon
as $\varphi$ is inner (w.l.o.g. $\varphi(0)=0$ again) and $h \in
H^2$ satisfies $\langle h,h\varphi^n \rangle=0$ for all $n \ge 1$
(which is a weaker condition than $h \in K_\varphi$); for
\[
\left\| \frac{h}{1-\overline w \varphi} \right \|_2^2= \left\|
\sum_{k=0}^\infty h\overline w^k \varphi^k \right\|_2^2 = \|h\|^2
\sum_{k=0}^\infty |w|^{2k},
\]
since $\langle h\varphi^m,h\varphi^n \rangle = 0$ for $m \ne n$.
Then by
 (\ref{eq:zen}) it follows that $W_{h,\varphi}$ is bounded.

In order to make a more systematic analysis of the weights that
induce a bounded weighted composition operator, we make the
following definition.
\begin{definition}
For $\varphi: \DD \to \DD$ analytic, the {\em multiplier space\/}
of $\varphi$ is defined by
\[
\mathcal{M}_\varphi = \{ h \in \mathcal{H}^2: \,
W_{h,\varphi}:=T_h C_\varphi \hbox{ is bounded} \}.
\]
\end{definition}

Evidently, $\mathcal{H}^\infty \subseteq \mathcal{M}_\varphi
\subseteq \mathcal{H}^2$ for all analytic self-maps $\varphi$ of
the unit disc. It is easily verified that $\mathcal{M}_\varphi$ is
a Banach space with the norm
$\|h\|_{\mathcal{M}_\varphi}=\|W_{h,\varphi}\|$.

Contreras and Hern\'andez-D\'{\i}az
\cite{CH1} and Matache \cite{matache} showed that
$\mathcal{M}_\varphi=\mathcal{H}^{\infty}$ if and only if
$\varphi$ is a finite Blaschke product.
In Corollary~\ref{cor:2.8} we shall present a stronger result.

The other extreme situation is the following.

\begin{theorem}
$\mathcal{M}_{\varphi}=\mathcal{H}^2$ if and only if
$\|\varphi\|_\infty < 1$.
\end{theorem}

\beginpf
If $\|\varphi\|_\infty < 1$, and $\|f\|_2=1$, then
\begin{eqnarray*}
\|C_\varphi f\|_\infty &=&  \left\|\sum_{n=0}^\infty \hat f(n)\varphi^n \right\|_\infty \\
& \le & \left(\sum_{n=0}^\infty |\hat f(n)|^2 \right)^{1/2}\left(
\sum_{n=0}^\infty \|\varphi\|_\infty^{2n} \right)^{1/2},
\end{eqnarray*}
and we conclude easily that $C_\varphi$ maps $\mathcal{H}^2$
boundedly into $\mathcal{H}^{\infty}$, whence
$\mathcal{M}_{\varphi}=\mathcal{H}^2$.

Conversely, if $\mathcal{M}_{\varphi}=\mathcal{H}^2$, then clearly
$C_\varphi$ maps $\mathcal{H}^2$  into $\mathcal{H}^{\infty}$, and
the mapping is bounded by the closed graph theorem. In particular
the functions $$\frac{\ds (1-|a|^2)^{1/2}}{\ds 1- \overline a
\varphi(z)},\qquad (z\in \mathbb{D}),$$ are uniformly bounded for
$a \in \DD$. However, if $|\varphi(z_0)|=1-\delta$ for some $z_0
\in \DD$ and $0< \delta < 1$, then for a suitable choice of $a$ of
modulus $1-\delta$ we obtain
\[
\left \| \frac{\ds (1-|a|^2)^{1/2}}{\ds 1- \overline a \varphi}
\right\|_\infty \ge \frac{\ds (1-(1-\delta)^2)^{1/2}}{\ds 1-
(1-\delta)^2} = (2\delta - \delta^2)^{-1/2},
\]
which tends to $\infty$ as $\delta \to 0$. Hence if
$\mathcal{M}_{\varphi}=\mathcal{H}^2$ we must have
$\|\varphi\|_\infty < 1$.
\endpf

There is a close relation between the angular derivative of
$\varphi$ and the boundedness of $W_{h,\varphi}$, as we explain in the next section.

\subsection{Relation with the angular derivative}

Once more, let $k_{w}$ be the reproducing kernel at $w$ in
$\mathcal{H}^2$, that is,
$$
k_{w}(z)=\frac{1}{1-\overline{w}\, z}, \qquad (z\in \mathbb{D}).
$$
If $W^{\ast}_{h,\varphi}$ denotes the adjoint of $W_{h,\varphi}$,
it is not difficult to see  that
$$
W_{h,\varphi}^{\ast}\, k_{w}= \overline{h(w)}\;  k_{\varphi(w)},
\qquad (w\in \mathbb{D}).
$$

Hence, if $W_{h,\varphi}$ is bounded, it follows that
\begin{equation}\label{eq1}
|h(w)|\leq \|W_{h,\varphi}\|\; \left (
\frac{1-|\varphi(w)|^2}{1-|w|^2} \right )^{1/2}
\end{equation}
for any $w\in \mathbb{D}$. Note that from equation (\ref{eq1}) one
gets another proof of the fact that
if $\varphi$ is a finite Blaschke product then ${\mathcal M}_\varphi=\mathcal{H}^{\infty}$;
for then we have
$$
\sup_{w\in \mathbb{D}}\frac{1-|\varphi(w)|^2}{1-|w|^2}<\infty,
$$
(as noted in \cite{matache}), and so  $h\in
\mathcal{H}^{\infty}$. The next result generalizes this observation
giving a necessary condition for boundedness and compactness of
weighted composition operators.

\begin{proposicion} \label{prop:angularderivative}
Let $\varphi$ be an analytic self-map of $\mathbb{D}$. Let
$E_{\varphi}=\{\zeta\in \mathbb{T}:\; \varphi \mbox{ has finite
angular derivative at } \zeta\}$. Then:
\begin{enumerate}
\item[i)] If $W_{h,\varphi}$ is bounded on $\mathcal{H}^p$ for some
$1\leq p <\infty$, then $h$ is pointwise bounded on every Stolz domain whose vertex is a point of $E_{\varphi}$.

\item[ii)]If $W_{h,\varphi}$ is compact on $\mathcal{H}^p$ for some
$1\leq p <\infty$, then $h$ tends to zero as $|z| \to 1$ in any Stolz domain whose vertex is a point  of
$E_{\varphi}$.
\end{enumerate}
\end{proposicion}

\beginpf
Note that if $\varphi$ has a finite angular derivative at
$\zeta\in \mathbb{T}$, then by the Julia--Carath\'eodory Theorem, the
nontangential limit satisfies
$$
\angle\lim_{ w\to \zeta}\frac{1-|\varphi(w)|}{1-|w|}<\infty
$$
(see \cite[pp. 51]{CM95} or \cite{ShBook}). Indeed, it holds that
$$
|\varphi^{\prime}(\zeta)|=\lim_{r\to 1}|\varphi^{\prime}(r\,
\zeta)|= \angle\lim_{ w\to \zeta}\frac{1-|\varphi(w)|}{1-|w|}=
\liminf_{ w\to \zeta}\frac{1-|\varphi(w)|}{1-|w|}.
$$
This along with (\ref{eq1}) and a standard argument on
reproducing kernels yield the statement of the proposition in
$\mathcal{H}^2$. The result in $\mathcal{H}^p$ is analogous.
\endpf

At this point, one may ask whether Proposition
\ref{prop:angularderivative} is also sufficient for boundedness of
weighted composition operators. The next example shows that this is no
longer true.

\begin{example}
{\rm Let $\varphi$ be the map $\varphi(z) = 1-\sqrt{1-z}$. Note
that $\varphi$ sends $\DD$  onto a subdomain of $\DD$ shaped
roughly like a (two dimensional) ice cream cone with vertex at 1.
It is easy to see that this map has finite angular derivative
nowhere on the unit circle. So for any $h\in \mathcal{H}^2$ the
boundedness requirement on $h$ given by Proposition
\ref{prop:angularderivative}(i) holds trivially.

Let us consider the $\mathcal{H}^2$ function $h(z) = (1-z)^{-3/8}$
and the weighted composition operator operator $W_{h,\varphi}$.
Observe that if 
$
f(z)=(1-z)^{-1/4}$, then
$$
(W_{h,\varphi}f)(z)=(1-z)^{-1/2},
$$
so $W_{h,\varphi}$ does not take $\mathcal{H}^2$ into itself. Thus
$W_{h,\varphi}$ satisfies the conclusion of Proposition
\ref{prop:angularderivative}(i), but not the hypothesis.} {\par
\hfill \rule{2.2mm}{2.2mm}\break}
\end{example}

If $\varphi$ is an inner function, there is a well known
characterization of the set where $\varphi$ has finite angular
derivative, essentially due to Carath\'eodory. Let
\begin{eqnarray}\label{eq:decomposition-canonique}
\varphi (z)=e^{i\alpha}z^N\prod_{n\geq 1}\frac{|a_n|}{a_n}
\frac{a_n -z}{1-\overline{a_n}z}\exp\left(-\int_{\mathbb T}\frac
{\zeta +z}{\zeta-z}\,d\mu(\zeta)\right),
\end{eqnarray}
where $a \in \RR$ and $N \in \NN \cup \{0\}$, be its canonical
representation in terms of the zeroes $\{a_n\}_n$ and a singular
measure $\mu$. We define the {\em Ahern--Clark set\/} $E_\varphi$
\cite{ahern70} by:
\[
E_\varphi:=\left\{\zeta\in\TT\,:\; \sum_{n\geq 1}
\frac{1-|a_n|^2}{|\zeta-a_n|^2}+2\int_{\TT}
\frac{d\mu(t)}{|t-\zeta|^2} <+\infty \right\}.
\]
\begin{proposition} (See, for example, \cite{CFP}.)
Let $\varphi$ be an inner function and $\zeta_0\in\TT$. Then the
following
assertions are equivalent:\\
(i) $\varphi$ has an angular derivative in the sense of
Carath\'eodory at $\zeta_0$.\\
(iii) $\zeta_0\in E_\varphi$.
\end{proposition}

It is known that $E_\varphi \supseteq \TT \setminus
\sigma(\varphi)$, where $\sigma(\varphi)$ denotes
the {\em spectrum\/} of $\varphi$, namely, $\sigma(\varphi)=\clos \{a_n\}_n \cup \supp \mu$. In general we do not have equality.\\

The following result may be seen as a generalization of the result
that $\mathcal{M}_{\varphi}=\mathcal{H}^{\infty}$ whenever
$\varphi$ is a finite Blaschke product (in which case
$E_\varphi=\TT$ and $\sigma(\varphi)=\emptyset$).

\begin{corollary}\label{cor:2.8}
Suppose that $\varphi$ is an inner function whose zero set, if
infinite, forms a finite  union of subsequences tending
nontangentially to points of $\TT$, and whose singular part is
given by a measure of finite support. Then any function $h \in
\mathcal{M}_{\varphi}$ is essentially bounded on all relatively
compact subsets of $E_\varphi=\TT \setminus \sigma(\varphi)$.
\end{corollary}

\beginpf
The fact that $E_\varphi=\TT \setminus \sigma(\varphi)$ under the
given hypotheses is shown in \cite[Thm.~5.3]{CFP}. On sufficiently
small neighbourhoods of any relatively compact subset $K$ of $\TT
\setminus \sigma(\varphi)$, one may extend $\varphi$ analytically
across $\TT$; we have $|\varphi(\zeta)|=1$ and $|\varphi'(\zeta)|$
is uniformly bounded for $\zeta\in K$. It is now easily verified
that $(1-|\varphi(w)|^2)/(1-|w|^2)$ is uniformly bounded as $w$
approaches $K$. The result follows from (\ref{eq1}).
\endpf

\begin{example}{\rm
The following example shows that the condition of Corollary \ref{cor:2.8}
is not sufficient for membership of $\mathcal{M}_{\varphi}$, even for $h \in \mathcal{H}^2$. Let
\[
\varphi(z)= \exp \left( -\frac{1+z}{1-z}\right), \qquad (z \in \DD),
\]
which is a singular inner function with $\sigma(\varphi)=\{1\}$. For $w=re^{i\alpha} \in \DD$, it is straightforward to verify that
\[
\frac{1-|\varphi(w)|^2}{1-|w|^2} =
\frac{1}{1-r^2} \left( 1 - \exp\left(-2\frac{1-r^2}{|1-re^{i\alpha}|^2} \right) \right),
\]
which, as $r \to 1$, tends to $2 /|\alpha|^2$.
Hence by (\ref{eq1}),  for each $h \in \mathcal{M}_{\varphi}$ we have that $|h(e^{i\alpha})|=O(|\alpha|^{-1})$ as $\alpha \to 0$.
It is easy to construct examples of functions $h \in {\mathcal H}^2$, bounded on any closed subarc of $\TT$
disjoint from $\{1\}$, that nonetheless do not satisfy this condition, and hence induce unbounded
operators $W_{h,\varphi}$.
}
\end{example}

\subsection{Schatten class operators}

In this section we provide stronger versions of the two main theorems of
\cite{gajath}. To begin with, in \cite[Thm.~2]{gajath}, it was shown that if
$\|\varphi\|_\infty < 1$, then
for every $h \in
\mathcal{H}^2$ the weighted composition operator $T=W_{h,\varphi}$
is compact. Indeed, much more is true, as the next result shows.

\begin{theorem}
Suppose that $\|\varphi\|_\infty < 1$. Then for every $h \in
\mathcal{H}^2$ the weighted composition operator $T=W_{h,\varphi}$
is trace-class (i.e., lies in $S_1$).
\end{theorem}
\beginpf
Note that for $w \in \DD$
\[
 \|T \widetilde{\dot k_w}\|_2 \approx (1-|w|^2)^{3/2}\left \|\frac{h\varphi}{(1-\overline
 w\varphi)^2}\right\|_2
\]
and since $|1-\overline w\varphi| \ge 1-\|\varphi\|_\infty$, this
is bounded by a constant multiple of $(1-|w|^2)^{3/2}\|h\|_2$. Now
a result of Harper and Smith \cite[Thm.~3.1]{zenmartin} asserts that
for $1 \le p < 2$, a  sufficient condition for a
weighted composition operator $T$ to lie in $S_p(\mathcal{H}^2)$
is
\begin{equation}\label{eq:zmintegral}
\int_{\DD} \|T \widetilde{\dot k_w} \|^p \frac{dA(w)}{(1-|w|^2)^2}
< \infty,
\end{equation}
where $\widetilde{ \dot k_w}=\dot k_w/\|\dot k_w\|$ and $\dot
k_w(z)=\displaystyle\frac{z}{(1-\overline w z)^2}$, so that
$\|\dot k_w\| \approx (1-|w|^2)^{-3/2}$.
We see now that the integral in (\ref{eq:zmintegral}) converges, even for $p=1$, and the result follows.
\endpf

The other main result of \cite{gajath}, its Theorem 1, asserts
that if $W_{h,\varphi}$ is bounded on $\mathcal{H}^2$, if
$\varphi$ lies in the disc algebra, and if $h=0$ on $\{\zeta\in
\TT: |\varphi(\zeta)|=1 \}$ and  $h$ is continuous at all such
points $\zeta$, then $W_{h,\varphi}$ is compact.

Here again, we may give a significantly stronger result.

\begin{theorem}
Suppose that $\varphi: \DD \to \DD$ is holomorphic and that  for
some $\delta>0$ and $c_\delta>0$ the function $h \in
\mathcal{H}^2$ satisfies $|h(z)| \le c_\delta$ a.e.\ on the set \[
 A_\delta:=\{z \in \TT: \, |\varphi(z)| \ge 1-\delta\}.
\]
 Then $W_{h,\varphi}$ is a bounded operator.
If, in addition, $c_\delta\to 0$ as $\delta \to 0$, then
$W_{h,\varphi}$ is a compact operator.
\end{theorem}
\beginpf
We estimate the $\mathcal{H}^2$ norm of
$(1-|w|^2)^{1/2}h/(1-\overline w\varphi)$ for $w \in \DD$.  Now
\[
\int_{A_\delta} \frac{(1-|w|^2)|h|^2}{|1-\overline w \varphi|^2}
\le c_\delta^2(1-|w|^2) \|C_\varphi k_w\|_2^2 \le c_\delta^2
\|C_\varphi\|^2,
\]
whereas
\[
\int_{\TT\setminus A_\delta}  \frac{(1-|w|^2)|h|^2}{|1-\overline w
\varphi|^2} \le \frac{(1-|w|^2) \|h\|^2_2 }{1-|w|(1-\delta)^2} \le
\frac{\|h\|_2^2 (1-|w|^2)}{\delta^2}.
\]
Thus
\[
\left\| \frac{(1-|w|^2)^{1/2}h}{1-\overline w\varphi}\right\|_2^2
\le c_\delta^2 \|C_\varphi\|^2 + \|h\|_2^2 (1-|w|^2)/\delta^2.
\]

This is uniformly bounded, independently of $w$, and so
$W_{h,\varphi}$ is a bounded operator; moreover, if $c_\delta\to
0$ as $\delta \to 0$, then given $\epsilon>0$, we may choose
$\delta>0$ sufficiently small
 such that
 $2c_\delta^2 \|C_\varphi\|^2<\epsilon$, and then find $\eta>0$ such that, if  $|w|>1-\eta$, one has $\|h\|_2^2 (1-|w|^2)/\delta^2 < c_\delta^2 \|C_\varphi\|^2$.
Hence condition (\ref{eq:zencompact}) is satisfied, and
$W_{h,\varphi}$ is compact.
\endpf

\section{Weighted Bergman spaces}
\label{sec:3}

For $\alpha>-1$ the weighted Bergman space
$\mathcal{A}^2_\alpha=\mathcal{A}^2_\alpha(\DD)$ consists of all
analytic functions in $\DD$ for which the norm, given by
\[
\|f\|_{\mathcal{A}^2_\alpha} =\left( \frac{1}{\pi} \int_\DD
|f(z)|^2 (1-|z|^2)^\alpha \, dA(z) \right)^{1/2}
\]
is finite. It is a reproducing kernel Hilbert space with kernel
function given by
\[
h_w(z)=\frac{\alpha+1}{(1-\overline wz)^{\alpha+2}}
\]
 (cf. \cite[p.~27]{CM95} and \cite[p.~135]{zhu}).
We shall also require the normalized functions $\widetilde
h_w=h_w/\|h_w\|_{\mathcal{A}^2_\alpha}$, given by
\[
\widetilde h_w(z)= \frac{(1-|w|^2)^{(\alpha+2)/2}}{(1-\overline
wz)^{\alpha+2}}.
\]

A necessary and sufficient condition for boundedness of
$W_{h,\varphi}$ on the Bergman space
$\mathcal{A}^2(\DD)=\mathcal{A}^2_0(\DD)$ is given in
\cite[Thm.~3.1]{KP}, namely that the measure $\nu_{h,\varphi}$
defined by
\[
\nu_{h,\varphi}(E)= \int_{\varphi^{-1}(E)} |h(z)|^2 \, dA(z),
\]
should be a Hastings--Carleson measure, i.e., the embedding $J:
\mathcal{A}^2(\DD) \to L^2(\DD,d\nu_{h,\varphi})$ should be
bounded. Such measures were characterised in
\cite{hastings,luecking}. These can also be tested on normalized
reproducing kernels $\widetilde h_w$, for $w \in \DD$ as in
\cite[Thm.~7.5]{zhu}.

Zen Harper \cite[Thm.~3.3]{zen} shows (as part of a more general
theory, which we describe in Theorem~\ref{thm:bigzen}, below) that
$W_{h,\varphi}$ is bounded on $\mathcal{A}^2(\DD)$ if and only if
\beq\label{eq:zen2} \sup_{|w|<1} \left \|
\frac{(1-|w|^2)h}{(1-\overline w \varphi)^2} \right
\|_{\mathcal{A}^2(\mathbb{D})} < \infty. \eeq

This result generalizes to $\mathcal{A}^2_\alpha(\DD)$ as follows.

\begin{proposition}
Let $\varphi$ be a holomorphic self-map of the disc and $h \in
\mathcal{A}^2_\alpha(\DD)$. Then the following conditions are
equivalent:

(i) The weighted composition operator $W_{h,\varphi}$ is bounded
on $\mathcal{A}^2_\alpha(\mathbb{D})$.

(ii) The measure $\nu_{h,\varphi}$ defined by
\[
\nu_{h,\varphi}(E)= \int_{\varphi^{-1}(E)} |h(z)|^2
(1-|z|^2)^\alpha \, dA(z)
\]
is an $\mathcal{A}^2_\alpha$-Carleson measure, in the sense that
the canonical embedding $J: \mathcal{A}^2_\alpha(\DD) \to
L^2(\DD,d\nu_{h,\varphi})$ is bounded.

(iii) One has \beq\label{eq:postzen} \sup_{|w|<1} \left \|
\frac{(1-|w|^2)^{1+\alpha/2}h}{(1-\overline w \varphi)^{\alpha+2}}
\right \|_{\mathcal{A}^2_\alpha} < \infty. \eeq

Moreover, $W_{h,\varphi}$ is compact if and only if
$\nu_{h,\varphi}$ is a vanishing Carleson measure, or,
equivalently, if \beq\label{eq:postzen2}
 \left \| \frac{(1-|w|^2)^{1+\alpha/2}h}{(1-\overline w \varphi)^{\alpha+2}} \right \|_{\mathcal{A}^2_\alpha} \to 0 \qquad \hbox{as} \quad |w| \to 1.
\eeq
\end{proposition}

\beginpf
The equivalence of (i) and (ii) follows from the formula
\[
\int_\DD |h(z)|^2 |f(\varphi(z))|^2 \, (1-|z|^2)^\alpha \, dA(z) =
\int_\DD |f(z)|^2 \, d\nu_{h,\varphi}(z),
\]
adapting the calculation for $\mathcal{A}^2(\DD)$, given in
\cite{KP}, without significant alteration. Next, observing that
\[
\int_{\DD} |\widetilde h_w|^2 \, d\mu_{h,\varphi} =  \int_{\DD}
|h(z)|^2 |\widetilde h_w \circ \varphi(z)|^2\, (1-|z|^2)^{\alpha}
\, dA(z),
\]
and using the fact that $\mathcal{A}^2_\alpha$-Carleson measures
can be tested on normalized reproducing kernels (which may be
found in \cite[Sec.~7.2]{zhu}), we see that (iii) implies (ii).
Trivially, (i) implies (iii), this being the property  that the
weighted composition
operator is uniformly bounded on the normalized kernels.\\

The proof of the characterization of compactness is similar, and
is omitted.
\endpf

To conclude, we mention briefly how one may also derive  results
that provide characterizations of Schatten class weighted
composition operators on $\mathcal{A}^2_\alpha(\mathbb{D})$. The
key is the following result of Harper
\cite[Thm.~3.1,~Cor.~3.2]{zen}.

\begin{theorem}\label{thm:bigzen}
Let $T$ be a subnormal operator on a Hilbert space $\mathcal{H}$,
with spectral radius $r(T) \le 1$. Let $h \in \mathcal{H}$. Then
there exists a finite positive Borel measure $\mu$ on
$\overline\DD$ such that
\[
\|f(T)h\|^2 = \int_{\overline \DD} |f(z)|^2 \, d\mu(z)
\]
for all functions $f$ analytic in a neighbourhood of $\overline
\DD$.
\end{theorem}

The application to weighted composition operators is the
following. Let $T=T_\varphi$ denote the operator of multiplication
by $\varphi$ on $\mathcal{A}^2_\alpha(\mathbb{D})$, which is
subnormal, since $\mathcal{A}^2_\alpha(\mathbb{D})$ is a closed
subspace of $L^2(\DD,(1-|z|^2)^\alpha \, dA(z))$. Then $f(T)h=h(f
\circ \varphi) = W_{h,\varphi}f$. We thus have the following
corollary.

\begin{corollary}
For $\varphi: \DD \to \DD$ holomorphic and $h \in
\mathcal{A}^2_\alpha(\mathbb{D})$, there exists a finite positive
Borel measure $\mu$ on $\overline\DD$ such that
\[
\|W_{h,\varphi}f\|_{\mathcal{A}^2_\alpha}^2 = \int_{\overline \DD}
|f(z)|^2 \, d\mu(z)
\]
for all functions $f$ analytic in a neighbourhood of $\overline
\DD$. Thus $\|W_{h,\varphi}f\|=\|I_\mu f\|$, where $I_\mu:
\mathcal{A}_\alpha^2(\DD) \to L^2(\overline\DD,\mu)$ is the
canonical Carleson embedding. It follows that $W_{h,\varphi}$ is
respectively bounded, compact or of Schatten class $S_p$ if and
only if $I_\mu$ has the same property.
\end{corollary}

We shall omit the detailed calculations, which are entirely
analogous to those of \cite[Prop.~2.4]{zenmartin}, but it may now
be verified that the Schatten class membership of $I_\mu$, and
hence $W_{h,\varphi}$, can be tested on reproducing kernels.
Namely $I_\mu \in S_p(\mathcal{A}^2_\alpha(\mathbb{D}),L^2(\mu))$
for $1<p<\infty$ if and only if
\[
\int_\DD \|I_\mu\widetilde h_w\|^p \, \frac{dA(w)}{(1-|w|^2)^2} <
\infty,
\]
and $I_\mu \in S_1(\mathcal{A}^2_\alpha(\DD),L^2(\mu))$ for
$1<p<\infty$ if and only if
\[
\int_\DD \|I_\mu\widetilde {\dot h}_w\|^p \,
\frac{dA(w)}{(1-|w|^2)^2} < \infty.
\]
where we write ${\dot h}_w=\frac{\partial}{\partial\overline
w}h_w$ and $\widetilde {\dot h}_w= {\dot h}_w/\| {\dot h}_w\|$. To
within irrelevant constants, we have
\[
 {\dot h}_w(z)\approx \frac{z}{(1-\overline w z)^{\alpha+3}}
\qquad \hbox{and} \qquad \widetilde {\dot h}_w(z) \approx
\frac{(1-|w|^2)^{(\alpha+4)/2}z}{(1-\overline w z)^{\alpha+3}}.
\]

\end{document}